\documentclass[letterpaper,11pt]{article}
\usepackage[margin=1in]{geometry}  

\usepackage{bbm}
\usepackage{graphicx}
\usepackage{amsmath,amssymb,amsthm,amsfonts}

\usepackage{paralist}
\usepackage{bm}
\usepackage{xspace}
\usepackage{url}
\usepackage{prettyref}
\usepackage{boxedminipage}
\usepackage{wrapfig}
\usepackage{color}
\usepackage{xspace}

\usepackage{amsmath,amsthm,amsfonts,amssymb}
\usepackage{graphicx}

\usepackage{nicefrac}

\newtheorem*{definition*}{Definition}

\usepackage{subcaption}

\usepackage[utf8]{inputenc}

\newcommand{\diag}{\textsf{diag}}

\usepackage{xcolor}
\definecolor{expert}{HTML}{008000}
\definecolor{error}{HTML}{f96565}

\usepackage{color-edits}
\addauthor{sw}{blue}

\usepackage{thmtools}
\usepackage{thm-restate}

\usepackage{tikz}
\usetikzlibrary{arrows,calc} 
\newcommand{\tikzAngleOfLine}{\tikz@AngleOfLine}
\def\tikz@AngleOfLine(#1)(#2)#3{%
\pgfmathanglebetweenpoints{%
\pgfpointanchor{#1}{center}}{%
\pgfpointanchor{#2}{center}}
\pgfmathsetmacro{#3}{\pgfmathresult}%
}

\declaretheoremstyle[
    headfont=\normalfont\bfseries, 
    bodyfont = \normalfont\itshape]{mystyle} 

\usepackage{listings}
\usepackage{amsmath}
\usepackage{amsthm}
\usepackage{tikz}
\usepackage{caption}
\usepackage{mdwmath}
\usepackage{multirow}
\usepackage{mdwtab}
\usepackage{eqparbox}
\usepackage{multicol}
\usepackage{amsfonts}
\usepackage{tikz}
\usepackage{multirow,bigstrut,threeparttable}
\usepackage{amsthm}
\usepackage{bbm}
\usepackage{epstopdf}
\usepackage{mdwmath}
\usepackage{mdwtab}
\usepackage{eqparbox}
\usetikzlibrary{topaths,calc}
\usepackage{latexsym}
\usepackage{amssymb}
\usepackage{bm}
\usepackage{amssymb}
\usepackage{graphicx}
\usepackage{mathrsfs}
\usepackage{epsfig}
\usepackage{psfrag}
\usepackage[
            CJKbookmarks=true,
            bookmarksnumbered=true,
            bookmarksopen=true,
            colorlinks=true,
            citecolor=red,
            linkcolor=blue,
            anchorcolor=red,
            urlcolor=blue
            ]{hyperref}

\usepackage{comment}
\usepackage{mathtools}
\usepackage{blkarray}
\usepackage{multirow,bigdelim,dcolumn,booktabs}

\usepackage{xparse}
\usepackage{tikz}
\usetikzlibrary{calc}
\usetikzlibrary{decorations.pathreplacing,matrix,positioning}

\usepackage[T1]{fontenc}
\usepackage[utf8]{inputenc}
\usepackage{mathtools}
\usepackage{blkarray, bigstrut}
\usepackage{gauss}

\newcommand*{\BraceAmplitude}{0.4em}%
\newcommand*{\VerticalOffset}{0.5ex}%
\newcommand*{\HorizontalOffset}{0.0em}%
\newcommand*{\blocktextwid}{3.0cm}%
\NewDocumentCommand{\InsertLeftBrace}{%
	O{} 
	O{\HorizontalOffset,\VerticalOffset} 
	O{\blocktextwid} 
	m   
	m   
	m   
}{%
	\begin{tikzpicture}[overlay,remember picture]
	\coordinate (Brace Top)    at ($(#4.north) + (#2)$);
	\coordinate (Brace Bottom) at ($(#5.south) + (#2)$);
	\draw [decoration={brace, amplitude=\BraceAmplitude}, decorate, thick, draw=black, #1]
	(Brace Bottom) -- (Brace Top) 
	node [pos=0.5, anchor=east, align=left, text width=#3, color=black, xshift=\BraceAmplitude] {#6};
	\end{tikzpicture}%
}%
\NewDocumentCommand{\InsertRightBrace}{%
	O{} 
	O{\HorizontalOffset,\VerticalOffset} 
	O{\blocktextwid} 
	m   
	m   
	m   
}{%
	\begin{tikzpicture}[overlay,remember picture]
	\coordinate (Brace Top)    at ($(#4.north) + (#2)$);
	\coordinate (Brace Bottom) at ($(#5.south) + (#2)$);
	\draw [decoration={brace, amplitude=\BraceAmplitude}, decorate, thick, draw=black, #1]
	(Brace Top) -- (Brace Bottom) 
	node [pos=0.5, anchor=west, align=left, text width=#3, color=black, xshift=\BraceAmplitude] {#6};
	\end{tikzpicture}%
}%
\NewDocumentCommand{\InsertTopBrace}{%
	O{} 
	O{\HorizontalOffset,\VerticalOffset} 
	O{\blocktextwid} 
	m   
	m   
	m   
}{%
	\begin{tikzpicture}[overlay,remember picture]
	\coordinate (Brace Top)    at ($(#4.west) + (#2)$);
	\coordinate (Brace Bottom) at ($(#5.east) + (#2)$);
	\draw [decoration={brace, amplitude=\BraceAmplitude}, decorate, thick, draw=black, #1]
	(Brace Top) -- (Brace Bottom) 
	node [pos=0.5, anchor=south, align=left, text width=#3, color=black, xshift=\BraceAmplitude] {#6};
	\end{tikzpicture}%
}%

\usetikzlibrary{patterns}

\definecolor{cof}{RGB}{219,144,71}
\definecolor{pur}{RGB}{186,146,162}
\definecolor{greeo}{RGB}{91,173,69}
\definecolor{greet}{RGB}{52,111,72}

\theoremstyle{plain}
\newtheorem{thm}{Theorem}[section]

\newtheorem{lemma}[thm]{Lemma}

\newtheorem{corollary}[thm]{Corollary}

\def \bP {\mathbb{P}}
\def \bQ {\mathbb{Q}}
\def \bC {\mathbb{C}}
\def \bE {\mathbb{E}}
\def \bR {\mathbb{R}}

\def\1{\mathbbm{1}}

\usepackage{xspace}

\newcommand{\pth}[1]{\left( #1 \right)}
\newcommand{\qth}[1]{\left[ #1 \right]}
\newcommand{\sth}[1]{\left\{ #1 \right\}}
\newcommand{\bpth}[1]{\Bigg( #1 \Bigg)}
\newcommand{\bqth}[1]{\Bigg[ #1 \Bigg]}
\newcommand{\bsth}[1]{\Bigg\{ #1 \Bigg\}}

\newcommand{\Bern}{\text{\rm Bern}}

\newcommand{\Unif}{\text{\rm Unif}}

\definecolor{myblue}{rgb}{.8, .8, 1}
\definecolor{mathblue}{rgb}{0.2472, 0.24, 0.6} 
\definecolor{mathred}{rgb}{0.6, 0.24, 0.442893}
\definecolor{mathyellow}{rgb}{0.6, 0.547014, 0.24}

\newcommand{\calC}{{\mathcal{C}}}

\newcommand{\calE}{{\mathcal{E}}}

\newcommand{\calN}{{\mathcal{N}}}

\newcommand{\calP}{{\mathcal{P}}}

\newcommand{\calX}{{\mathcal{X}}}

\newcommand{\rmd}{\mathrm{d}}
\newcommand{\Perm}{\mathrm{Perm}}

\usepackage{cleveref}
\crefname{lemma}{Lemma}{Lemmas}
\Crefname{lemma}{Lemma}{Lemmas}
\crefname{thm}{Theorem}{Theorems}
\Crefname{thm}{Theorem}{Theorems}
\Crefname{assumption}{Assumption}{Assumptions}
\crefformat{equation}{(#2#1#3)}

\begin{document}
\title{Elementary Symmetric Polynomial Inequalities \\ for Centered Vectors and Matrices}
\author{Yanjun Han, Jonathan Niles-Weed\thanks{Courant Institute of Mathematical Sciences and Center for Data Science, New York University. Email: \url{yanjunhan@nyu.edu}, \url{jnw@cims.nyu.edu}}}
\maketitle

\begin{abstract}
We prove new inequalities for elementary symmetric polynomials (ESPs) for vectors that sum to zero, and for square matrices with zero row and column sums. We apply these results to obtain a unified upper bound on the mean-field approximation guarantee for permutation mixtures, as well as a sharp $\chi^2$ version of the de Finetti theorem for finite sequences over a small alphabet. The main proof ideas were developed by the GPT-5.5 Pro model.
\end{abstract}

\section{Introduction}
Given $n$ real or complex numbers $z_1,\dots,z_n$ and $0\le k\le n$, the elementary symmetric polynomial (ESP) of degree $k$ is defined as
\begin{align}\label{eq:ESP}
e_k(z) = \sum_{S\subseteq [n], |S|=k}\prod_{i\in S}z_i = \sum_{1\le i_1<\dots<i_k\le n} z_{i_1}\cdots z_{i_k}. 
\end{align}
When $z_1,\dots,z_n$ are real, a classical question is to understand the geometry of the set of \emph{attainable} tuples $\{(e_0(z),\dots,e_n(z)): z\in \bR^n\}$. For instance, if $n=2$, this set is $\{(e_0,e_1,e_2): e_0=1, 4e_2\le e_1^2\}$. To understand the geometry for general $n$, a large collection of inequalities was known. The most notable examples include \emph{Newton's inequality}
\begin{align*}
\frac{e_{k+1}(z)}{\binom{n}{k+1}}\frac{e_{k-1}(z)}{\binom{n}{k-1}} \le \frac{e_k(z)^2}{\binom{n}{k}^2}
\end{align*}
for real vectors, and \emph{Maclaurin's inequality}
\begin{align*}
\pth{\frac{e_k(z)}{\binom{n}{k}}}^{1/k} \le \pth{\frac{e_{k-1}(z)}{\binom{n}{k-1}}}^{1/(k-1)} \le \cdots \le \frac{e_1(z)}{n}
\end{align*}
for \emph{non-negative} real vectors. Replacing $z$ by $|z|$, these yield inequalities for complex vectors $z$ as well. 

The main goal of this paper is to control the magnitude of $e_k(z)$ for \emph{centered} vectors, namely vectors $z$ with $\sum_{i=1}^n z_i = 0$. By scaling, we will further assume that $\sum_{i=1}^n |z_i|^2 = n$, so that the average squared magnitude of the coordinates of $z$ is normalized to $1$. In this case, the classical consequence $|\frac{e_k(z)}{\binom{n}{k}}|^{1/k}\le (\frac{1}{n}\sum_{i=1}^n |z_i|^2)^{1/2}$ of Maclaurin's inequality gives that $|e_k(z)|\le \binom{n}{k}$, which is tight for \emph{non-centered} vectors such as $z={\bf 1}$. For centered vectors, it will turn out that this upper bound can be significantly improved to $O(\binom{n}{k}^{1/2})$, a scaling partially inspired by the randomized counterpart $\bE[|e_k(z)|^2] = \binom{n}{k}$ for all independent random variables $z_1,\dots,z_n$ with zero mean and unit variance. Such upper bounds also add to a recent line of work \cite{gopalan2014inequalities, meka2019pseudorandom, doron2020log, tao2023maclaurin} proving bounds on elementary symmetric polynomials given the values of two consecutive lower-order elementary symmetric polynomials.

The rest of this paper is organized as follows. We first present a new upper bound of $|e_k(z)|$ for all centered complex vectors $z$, and discuss how it improves upon existing bounds. We then extend the definition of $e_k$ to centered square matrices, namely matrices with zero row and column sums, and present how centering along both directions leads to a smaller magnitude of $e_k$ compared with non-centered and partially centered settings. Finally, we discuss applications of these inequalities to problems involving random permutation structures, including permutation mixtures and de Finetti-type theorems. The core ideas underlying the proofs of both inequalities were initially proposed by the GPT-5.5 Pro model. The authors refined the arguments, developed the applications, and wrote the paper.

\subsection{An inequality for centered vectors}
We first present an upper bound of $e_k(z)$ when $z$ is a centered complex vector. 
\begin{thm}\label{thm:ESP-vector}
Let $z\in \bC^n$ with $\sum_{i=1}^n z_i = 0$ and $\sum_{i=1}^n |z_i|^2 = n$, and $0\le k\le n$. Then
\begin{align*}
|e_k(z)| \le \sqrt{C\binom{n}{k}}
\end{align*}
for some universal constant $C>0$. An unoptimized constant is $C = 24$. 
\end{thm}

We give an overview of the existing literature on a collection of weaker results. For a centered real vector $z$, \cite[Theorem 2]{gopalan2014inequalities} gives an upper bound using Lagrangian multipliers: 
\begin{align*}
|e_k(z)|^2 \le \pth{6e\sqrt{\frac{n}{k}}}^{k} \le (6e)^k\sqrt{\binom{n}{k}}, 
\end{align*}
which enjoys the $\binom{n}{k}^{1/2}$ growth but has an additional exponential factor in $k$. This exponential factor was significantly reduced by a differentiation trick in \cite{tao2023maclaurin}, which we reproduce here for completeness. Let $f(w) = \prod_{i=1}^n (w-z_i) = \sum_{k=0}^n (-1)^{k}e_k(z)w^{n-k}$, then
\begin{align*}
\frac{f^{(n-k)}(w)}{n!/k!} = \sum_{\ell=0}^k (-1)^{\ell}\frac{\binom{k}{\ell}}{\binom{n}{\ell}}e_{\ell}(z)w^{k-\ell}. 
\end{align*}
Writing $w_1,\dots,w_k$ for the roots of $f^{(n-k)}(w)$, the LHS equals $\prod_{\ell=1}^k (w-w_{\ell})$. The Schoenberg conjecture (now the Malamud--Pereira theorem \cite{malamud2005inverse,pereira2003differentiators}) relates the roots of $f'$ to the centered roots of $f$, and we repeatedly apply it $(n-k)$ times to get 
\begin{align*}
\frac{1}{k}\sum_{\ell=1}^k |w_\ell|^2 \le \frac{k-1}{n-1}\cdot \frac{1}{n}\sum_{i=1}^n |z_i|^2 = \frac{k-1}{n-1}. 
\end{align*}
Finally, by AM-GM and the inequality $\binom{n}{k}\le (\frac{en}{k})^k$, it holds that
\begin{align*}
|e_k(z)| = \binom{n}{k} \prod_{\ell=1}^k |w_\ell| \le \binom{n}{k} \pth{\frac{1}{k}\sum_{\ell=1}^k |w_\ell|^2}^{k/2} \le \binom{n}{k} \pth{\frac{k-1}{n-1}}^{k/2} \le \sqrt{e^k \binom{n}{k}}. 
\end{align*}
As a result, this exponential factor reduces from $(6e)^k$ to $e^{k/2}$. It can be further removed via a saddle point analysis of a complex contour integral. Specifically, \cite[Lemma 2.2]{roos2015bobkov} (see also \cite[Part 1 of Theorem 4.3]{han2024approximate}) shows that for $1\le k\le n-1$, 
\begin{align*}
|e_k(z)| \le \pth{\frac{n^n}{k^k(n-k)^{n-k}}}^{1/2}\le C\pth{\frac{k(n-k)}{n}}^{1/4} \sqrt{\binom{n}{k}}, 
\end{align*}
which has no extra exponential factor in $k$ but a small polynomial factor of $k^{1/4}$. Finally, \cite[Part 2 of Theorem 4.3]{han2024approximate} shows that this polynomial factor can also be removed for real $z$: $|e_k(z)|\le \sqrt{10\binom{n}{k}}$. The proof establishes that the maximum of $|e_k(z)|$ over real $z$ is attained by a vector taking only two values using Lagrangian multipliers, and applies a more careful saddle point analysis for these special vectors. However, the first step critically relies on a property of real-rooted polynomials (cf. \cite[Fact B]{gopalan2014inequalities}), so such reduction to vectors with binary support is not immediate for complex $z$. \Cref{thm:ESP-vector} generalizes this result to all complex vectors (with a larger constant $C$), via a careful saddle point analysis with sublevel estimates. 

\subsection{An inequality for centered matrices}
We also consider an ESP-like quantity for matrices. For a square matrix $A\in \bC^{n\times n}$ and $1\le k\le n$, define
\begin{align}\label{eq:ESP-matrix}
e_k(A) = \frac{1}{k!}\sum_{\substack{S,T\subseteq [n] \\ |S|=|T|=k}} \Perm(A_{S,T}), 
\end{align}
where $\Perm(A) = \sum_{\pi\in S_n} \prod_{i=1}^n A_{i,\pi(i)}$ is the matrix permanent, and $A_{S,T}$ denotes the submatrix with row indices in $S$ and column indices in $T$. The above definition, including the normalizing constant $\frac{1}{k!}$, is related to ESPs of vectors via $e_k(A) = e_k(u)e_k(v)$ for a rank-$1$ matrix $A=uv^\top$. The following result presents an upper bound of $e_k(A)$ for a \emph{centered} matrix $A$. 

\begin{thm}\label{thm:ESP-matrix}
Let $A\in \bC^{n\times n}$ have zero row and column sums, i.e. $A{\bf 1} = A^\top{\bf 1}=0$, and be normalized such that $\|A\|_{\mathrm{F}} = (\sum_{i,j} |A_{ij}|^2)^{1/2} = n$. Then for $1\le k\le n$, 
\begin{align*}
|e_k(A)| \le B^k\binom{n}{k}
\end{align*}
for a universal constant $B>0$. An unoptimized choice is $B=8e^2$. 
\end{thm}

Similar to the vector case, the centered assumption is important: the all-ones matrix $J$ satisfies $\|J\|_{\mathrm{F}}=n$, but $e_k(J) = \binom{n}{k}^2$. In addition, it is important that the matrix $A$ be centered along \emph{both} directions of rows and columns. For example, if $A$ only has zero row sums, then the choice $A={\bf 1}z^\top$ for a centered vector $z$ would give $e_k(A) = \binom{n}{k}e_k(z)$, with $|e_k(z)|$ possibly of the order of $\binom{n}{k}^{1/2}$. 

Unlike \Cref{thm:ESP-vector} where we can even obtain the tight exponential factor $B=1$, in the matrix case it appears much more challenging to exploit both the row and column constraints and obtain a weaker form as shown in \Cref{thm:ESP-matrix}. Before \Cref{thm:ESP-matrix}, the only known approach is to apply \cite[Lemma 4.2]{han2024approximate} (i.e. Banach theorem for symmetric multilinear forms) to the row vectors of $A$ to reduce to the case where $A$ has identical rows, so that
\begin{align*}
|e_k(A)| \le C\binom{n}{k}^{3/2} \pth{\prod_{i=1}^n \frac{\| A_{i,\cdot} \|_2^2}{n} }^{1/2} \le C\binom{n}{k}^{3/2} \pth{\frac{\|A\|_{\mathrm{F}}^2}{n^2} }^{k/2} = C\binom{n}{k}^{3/2}. 
\end{align*}
However, this approach fundamentally ignores the property that all column sums are also zero. 

Another tempting, but ultimately unsuccessful, approach is to reduce \Cref{thm:ESP-matrix} to the vector inequality in \Cref{thm:ESP-vector} via tensorization. Define the symmetric $k$-linear form $\widetilde{\calE}_k: (\bC^n)^k\to \bC$ by
\begin{align*}
\widetilde{\calE}_k\pth{x^{(1)},\dots,x^{(k)}} := \frac{1}{k!}\sum_{\substack{i_1,\dots,i_k\in [n] \\ \text{all distinct}}} \prod_{\ell=1}^k x^{(\ell)}_{i_\ell},
\end{align*}
and let $\calE_k := \widetilde{\calE}_k|{H^k}$ denote its restriction to $H := \{z\in \bC^n: \sum_{i=1}^n z_i=0\} = {\bf 1}^\perp$. For every $z\in \bC^n$, we have $\widetilde{\calE}_k(z,\dots,z) = e_k(z)$. For the matrix case, define the \emph{Kronecker (or vertical) tensor product} of two multilinear forms as follows. Let $\Lambda: E^k\to \bC$ and $\Gamma: F^k\to \bC$ be $k$-linear forms on finite-dimensional complex Hilbert spaces $E$ and $F$. Then $\Lambda\boxtimes\Gamma: (E\otimes F)^k\to \bC$ is the $k$-linear form determined on pure tensors by
\begin{align*}
(\Lambda\boxtimes\Gamma)(x_1\otimes y_1,\dots,x_k\otimes y_k) := \Lambda(x_1,\dots,x_k)\,\Gamma(y_1,\dots,y_k).
\end{align*}
Identify $\bC^{n\times n}$ isometrically with the tensor product $\bC^n\otimes \bC^n$ by associating each $A = (A_{ij})$ with the tensor $
a_A := \sum_{i,j=1}^n A_{ij}\, e_i\otimes e_j$, then $\|a_A\|_2 = \|A\|_{\mathrm{F}}$, and $a_A\in H\otimes H$ whenever $A{\bf 1}=A^\top {\bf 1}= 0$. It can then be verified that
\begin{align*}
e_k(A) = \calE_k \boxtimes \calE_k(a_A, \dots, a_A). 
\end{align*}
As a result, a stronger statement $|e_k(A)|\le 24\binom{n}{k}$ in \Cref{thm:ESP-matrix} would follow from \Cref{thm:ESP-vector} \textit{if} the spectral norm were multiplicative under Kronecker tensor products, i.e.
$
\|\Lambda\boxtimes\Gamma\| = \|\Lambda\|\, \|\Gamma\|$, with $\|\Lambda\| := \sup_{\|x_1\|_2=\dots=\|x_k\|_2=1} |\Lambda(x_1,\dots,x_k)|$. Unfortunately, this tempting claim is false.\footnote{Though multiplicativity of the spectral norm has been claimed in the literature~\cite[Proposition 3.3]{derksen2016nuclear}, the proof there has a gap.} Indeed, for $\Delta_k = \det[x_1, \dots, x_k]$ on $(\bC^k)^k$, it holds that $\|\Delta_k\|=1$ but $\|\Delta_k\boxtimes \Delta_k\|\ge \frac{k!}{k^{k/2}}$, as witnessed by evaluating the form at $(\Omega_k,\dots,\Omega_k)$, with $\Omega_k = \frac{1}{\sqrt{k}}\sum_{i=1}^k e_i\otimes e_i$.

Our proof of \Cref{thm:ESP-matrix} will rely on a carefully chosen generating function that respects the centering assumption along both directions.

\subsection{Applications}
\paragraph{Spectral upper bounds of matrix permanent.} Both inequalities find applications in the permanent upper bound for PSD and doubly stochastic matrices. Let $M\in \bR^{n\times n}$ be a PSD and doubly stochastic matrix, which must admit eigenvalues $1=\lambda_1\ge \lambda_2 \ge \cdots \ge \lambda_n\ge 0$. The well-known van der Waerden inequality \cite{van1926aufgabe,egorychev1981solution,falikman1981proof} states that $\Perm(M)\ge \frac{n!}{n^n}$, with equality iff $M = \frac{J}{n}$ is a scaled all-ones matrix (or equivalently, when $\lambda_1 = 1$ and $\lambda_2=\cdots=\lambda_n=0$). The next corollary presents two different upper bounds of $\Perm(M)$, from \Cref{thm:ESP-vector} and \ref{thm:ESP-matrix} respectively, in terms of the non-leading eigenvalues $\lambda_2,\dots,\lambda_n$. 

\begin{corollary}\label{cor:permanent}
Under the above setting, it holds that
\begin{align*}
\frac{n^n}{n!}\Perm(M) - 1 \le \sum_{k=2}^n \min\bsth{ C\sum_{ \substack{k_2,\dots,k_n\ge 0 \\ k_2+\dots+k_n=k }} \lambda_2^{k_2}\cdots\lambda_n^{k_n}, B^k\pth{\lambda_2^2+\cdots+\lambda_n^2}^{k/2} },
\end{align*}
where the universal constants $C,B$ are given in \Cref{thm:ESP-vector} and \ref{thm:ESP-matrix}, respectively. 
\end{corollary}

We discuss how \Cref{cor:permanent} improves over existing results in \cite{han2024approximate}. The first inequality
\begin{align*}
   \frac{n^n}{n!}\Perm(M) - 1 \le C \sum_{k=2}^n \sum_{ \substack{k_2,\dots,k_n\ge 0 \\ k_2+\dots+k_n=k }} \lambda_2^{k_2}\cdots\lambda_n^{k_n} 
\end{align*}
unifies and improves all three upper bounds of \cite{han2024approximate} (ignoring differences in the constant $C$): 
\begin{enumerate}
    \item The first upper bound of \cite{han2024approximate} reads
    \begin{align}\label{eq:known-UB-1}
    \frac{n^n}{n!}\Perm(M) - 1 \le C\sum_{k=2}^n (\lambda_2+\cdots+\lambda_n)^k, 
    \end{align}
    which was established using the counterpart of \Cref{thm:ESP-vector} for \emph{real} vectors $z\in \bR^n$. Since $\sum_{ \substack{k_2,\dots,k_n\ge 0 \\ k_2+\dots+k_n=k }} \lambda_2^{k_2}\cdots\lambda_n^{k_n}\le (\lambda_2+\cdots+\lambda_n)^k$, this bound is never better than \Cref{cor:permanent}.
    \item The second upper bound of \cite{han2024approximate} reads
    \begin{align*}
    \frac{n^n}{n!}\Perm(M) - 1 \le \sum_{k=1}^n \sum_{ \substack{k_2,\dots,k_n\ge 0 \\ k_2+\dots+k_n=k }} \lambda_2^{k_2}\cdots\lambda_n^{k_n}, 
    \end{align*}
    where the sum critically starts from $k=1$. In the local regime where $\lambda_2,\dots,\lambda_n$ are small, this upper bound exhibits a linear dependence on $(\lambda_2,\dots,\lambda_n)$ rather than quadratic.
    \item The third upper bound of \cite{han2024approximate} refines the second by starting from $k=2$, but with an additional polynomial factor: 
    \begin{align*}
    \frac{n^n}{n!}\Perm(M) - 1 \le C\sum_{k=2}^n \sqrt{k}\sum_{ \substack{k_2,\dots,k_n\ge 0 \\ k_2+\dots+k_n=k }} \lambda_2^{k_2}\cdots\lambda_n^{k_n}.  
    \end{align*}
    \Cref{cor:permanent} removes the extra $\sqrt{k}$ factor. 
\end{enumerate}
These upper bounds were established via three related but different approaches in \cite{han2024approximate}, and in general they are incomparable to each other. Therefore, our new bound in \Cref{cor:permanent} unifies and improves all results of \cite{han2024approximate}. 

\paragraph{Permutation mixture.} In \cite{han2024approximate}, the permanent upper bounds were applied to understanding the mean-field approximation of \emph{permutation mixtures}. Specifically, let $P_1,\dots,P_n\in \calP$ be $n$ probability distributions over the same space. The permutation mixture, as well as its i.i.d. approximation, is defined as
\begin{align}\label{eq:permutation-mixture}
\bP_n = \bE_{\pi\sim \Unif(S_n)}\bqth{ \bigotimes_{i=1}^n P_{\pi(i)} }, \quad \bQ_n = \bpth{\frac{1}{n}\sum_{i=1}^n P_i}^{\otimes n}. 
\end{align}
In words, the permutation mixture $\bP_n$ is the joint distribution of $n$ random draws from a noisy population $X_1\sim P_1,\dots,X_n\sim P_n$ without replacement, and $\bQ_n$ is the counterpart with replacement. The central result of \cite{han2024approximate} is a (usually dimension-free) upper bound on the statistical distance between $\bP_n$ and $\bQ_n$: 
\begin{align*}
\chi^2(\bP_n \| \bQ_n) \le \min\bsth{ C\frac{(\lambda_2(M)+\cdots+\lambda_n(M))^2}{(1-\lambda_2(M)-\cdots-\lambda_n(M))_+}, \prod_{i=2}^n \frac{1}{1-\lambda_i(M)}-1 }, 
\end{align*}
where $M=M(P_1,\dots,P_n)\in \bR^{n\times n}$ is a channel overlap matrix generated by $P_1,\dots,P_n$ (cf. \cite[Lemma 5.1]{han2024approximate} or \eqref{eq:channel-overlap} for its definition) which is doubly stochastic and PSD. Here the first bound illustrates the optimal quadratic dependence on $\sum_{i=2}^n \lambda_i(M)$ when it is small, but becomes vacuous when it is large. The second bound remains meaningful even when $\lambda_i(M)$ is close to $1$, but locally it suggests a linear dependence on $\sum_{i=2}^n \lambda_i(M)$ rather than quadratic. \Cref{cor:permanent} gives a ``best-of-both-worlds'' upper bound that improves over both results. 

\begin{corollary}\label{cor:chi-squared}
Under the above setting, it holds that
\begin{align*}
\chi^2(\bP_n \| \bQ_n) \le C\bpth{ \prod_{i=2}^n \frac{1}{1-\lambda_i(M)} - 1 - \sum_{i=2}^n \lambda_i(M) }. 
\end{align*}
Here the universal constant $C$ is given in \Cref{thm:ESP-vector}. 
\end{corollary}

\paragraph{Finite de Finetti theorem.} Instead of comparing the full joint distributions $\bP_n$ and $\bQ_n$, one may also compare their $k$-dimensional marginals, denoted by $\bP_{k,n}$ and $\bQ_{k,n}$, respectively. This question, commonly studied under the framework of the \emph{finite de Finetti theorem}, has inspired a rich line of classical and recent work in probability and information theory \cite{de1929funzione,diaconis1980finite,stam1978distance,gavalakis2021information,johnson2024relative}. A natural objective is to identify conditions on $k$ under which the statistical distance between $\bP_{k,n}$ and $\bQ_{k,n}$ is $o(1)$; equivalently, when sampling without replacement is statistically indistinguishable from sampling with replacement. For Dirac measures $P_1,\dots,P_n$, \cite{diaconis1980finite} shows that $k=o(\sqrt{n})$ is sufficient, while \cite{stam1978distance} improves this condition to $k=o(n/\sqrt{|\calX|})$ when $P_1,\dots,P_n$ are supported on a small alphabet $\calX$. By contrast, \cite{han2024approximate} shows that, in the noisy setting where no two measures among $P_1,\dots,P_n$ are mutually singular, even $k=o(n)$ may be sufficient.

The second inequality of \Cref{cor:permanent}, namely
\begin{align*}
\frac{n^n}{n!}\Perm(M) - 1 \le \sum_{k=2}^n B^k\pth{\lambda_2^2 + \cdots + \lambda_n^2}^{k/2},
\end{align*}
improves upon the existing upper bound \eqref{eq:known-UB-1} and provides a unified treatment of the above scenarios.  

\begin{corollary}\label{cor:deFinetti}
Let $M=M(P_1,\dots,P_n)\in \bR^{n\times n}$ be the channel overlap matrix generated by $P_1,\dots,P_n$, and ${\sf T} := \sum_{i=2}^n \lambda_i(M)$. There exists a universal constant $C>0$ such that
\begin{align*}
\chi^2(\bP_{n,k} \| \bQ_{n,k}) \le \frac{Ck^2}{n^2}\pth{ {\sf T}^2\wedge {\sf T} }, \quad \text{if } k \le \frac{n}{C\sqrt{{\sf T}^2\wedge {\sf T}}}. 
\end{align*}
In particular, a sufficient condition of $\chi^2(\bP_{n,k}\|\bQ_{n,k})=o(1)$ is $k=o(\frac{n}{\sqrt{{\sf T}}})$. 
\end{corollary}

To see how \Cref{cor:deFinetti} unifies the existing scenarios, note that ${\sf T}\le n-1$ always holds, since $\lambda_i\le 1$, and ${\sf T}\le |\calX|-1$ when $P_1,\dots,P_n$ are supported on a finite alphabet $\calX$; see \Cref{subsec:deFinetti} for a proof. We also note that \cite{stam1978distance,diaconis1980finite} primarily work with total variation (TV) distance or Kullback--Leibler (KL) divergence, whereas our result generalizes these bounds to the stronger $\chi^2$ divergence. In the noisy scenario, \cite[Theorem 1.4]{han2024approximate} shows that $k=o(\frac{n}{\sf T})$ is sufficient for $\chi^2(\bP_{n,k}\|\bQ_{n,k})=o(1)$, while \Cref{cor:deFinetti} improves the dependence to $k=o(\frac{n}{\sqrt{\sf T}})$ in the interesting regime ${\sf T}\ge 1$. In particular, both results show that if ${\sf T}=O(1)$, then $k=o(n)$ is sufficient.


\paragraph{Acknowledgement.} We thank Yuzhou Gu from OpenAI for using the GPT-5.5 Pro model to produce the proofs of \Cref{thm:ESP-vector} and \ref{thm:ESP-matrix}. We also thank Yuxiao Wen for helpful discussions during an early stage of \Cref{thm:ESP-matrix}.


\section{Proofs}
\subsection{Proof of \Cref{thm:ESP-vector}}
The cases $k=0$ and $k=n$ are trivial, so in the sequel we assume that $1\le k\le n-1$. We begin with the same contour integral in \cite[Proof of Theorem 4.3]{han2024approximate}: let $r=\sqrt{\frac{k}{n-k}}$, 
\begin{align*}
e_k(z) = \frac{1}{2\pi \mathrm{i}}\oint_{|w|=r}\frac{\prod_{i=1}^n(1+z_iw)}{w^{k+1}}\rmd w = \frac{1}{2\pi}\int_{0}^{2\pi} \frac{\prod_{i=1}^n (1+rz_ie^{\mathrm{i}\theta})}{r^ke^{\mathrm{i}k\theta}} \rmd \theta. 
\end{align*}
By triangle inequality, 
\begin{align}\label{eq:basic-ineq}
|e_k(z)| \le \frac{1}{2\pi r^k}\int_0^{2\pi} \prod_{i=1}^n |1+rz_ie^{\mathrm{i}\theta}| \rmd \theta. 
\end{align}
By the assumptions $\sum_{i=1}^n z_i=0$ and $\sum_{i=1}^n |z_i|^2 = n$, we have
\begin{align*}
\frac{1}{n}\sum_{i=1}^n |1+rz_ie^{\mathrm{i}\theta}|^2 = 1+r^2 + 2r\Re\pth{e^{\mathrm{i}\theta}\cdot \frac{1}{n}\sum_{i=1}^n z_i} = 1+r^2.
\end{align*}
Let
\begin{align*}
q(\theta) := \prod_{i=1}^n \frac{|1+rz_ie^{\mathrm{i}\theta}|}{\sqrt{1+r^2}}, 
\end{align*}
the AM-GM inequality applied to $q(\theta)^2$ yields $q(\theta)\le 1$ for all $\theta\in [0,2\pi]$. If we apply this uniform inequality to \eqref{eq:basic-ineq}, we would get the weaker upper bound $|e_k(z)|\le \frac{n^{n/2}}{k^{k/2}(n-k)^{(n-k)/2}}$ in \cite{roos2015bobkov,han2024approximate}. Therefore, it is important to provide a superlevel estimate of $q(\theta)$ for $\theta$ ranging over $[0,2\pi]$.  

To this end, we define two useful functions
\begin{align*}
A(\theta) := \sum_{i=1}^n \pth{\frac{|1+rz_ie^{\mathrm{i}\theta}|^2}{1+r^2} - 1}^2, \quad B(\theta) := \sum_{i=1}^n \pth{\frac{|1+rz_ie^{\mathrm{i}\theta}|}{\sqrt{1+r^2}} - 1}^2. 
\end{align*}
Expanding the square of $B(\theta)$, we have
\begin{align*}
B(\theta) = 2\sum_{i=1}^n \pth{1- \frac{|1+rz_ie^{\mathrm{i}\theta}|}{\sqrt{1+r^2}}}. 
\end{align*}
Therefore, if we apply AM-GM to $q(\theta)$ instead of $q(\theta)^2$, we get
\begin{align}\label{eq:q-upper-bound}
q(\theta)\le \pth{\frac{1}{n} \sum_{i=1}^n \frac{|1+rz_ie^{\mathrm{i}\theta}|}{\sqrt{1+r^2}} }^n \le \exp\pth{\sum_{i=1}^n \pth{\frac{|1+rz_ie^{\mathrm{i}\theta}|}{\sqrt{1+r^2}}}-1} = \exp\pth{-\frac{B(\theta)}{2}}. 
\end{align}
In addition, since
\begin{align*}
A(\theta) &=  \sum_{i=1}^n \pth{\frac{|1+rz_ie^{\mathrm{i}\theta}|}{\sqrt{1+r^2}} - 1}^2\pth{\frac{|1+rz_ie^{\mathrm{i}\theta}|}{\sqrt{1+r^2}} + 1}^2 \\
&\le  \sum_{i=1}^n \pth{\frac{|1+rz_ie^{\mathrm{i}\theta}|}{\sqrt{1+r^2}} - 1}^2(\sqrt{B(\theta)}+2)^2 = B(\theta)(\sqrt{B(\theta)}+2)^2, 
\end{align*}
we have $\sth{\theta: B(\theta)\le t} \subseteq \sth{\theta: A(\theta)\le t(\sqrt{t}+2)^2}$ for all $t\ge 0$. In other words, an upper bound on the lower tail of $A(\theta)$ controls the lower tail of $B(\theta)$ as well. 

The remaining step is to control the lower tail of $A(\theta)$. It is clear that $A(\theta)\ge 0$ is a trigonometric polynomial of degree $2$, i.e. $A(\theta)=a_0+\sum_{j=1}^2 (a_j\cos(j\theta)+b_j\sin(j\theta))$. In addition, 
\begin{align*}
\frac{1}{2\pi}\int_0^{2\pi} A(\theta)\rmd \theta &= \frac{1}{2\pi(1+r^2)^2}\sum_{i=1}^n \int_0^{2\pi} \pth{r^2(|z_i|^2-1)+2r\Re(z_ie^{\mathrm{i}\theta})}^2 \rmd \theta \\
&= \frac{1}{(1+r^2)^2}\sum_{i=1}^n\pth{r^4 (|z_i|^2-1)^2 + 2r^2|z_i|^2} \ge \frac{2nr^2}{(1+r^2)^2}. 
\end{align*}
We use the following version of Remez inequality, due to Ganzburg \cite{ganzburg2012remez}, to obtain a sublevel estimate of $A(\theta)$. 

\begin{lemma}\label{lemma:remez}
Let $T_d$ be a trigonometric polynomial of degree at most $d$ with real coefficients, and $E\subseteq [0,2\pi]$ be a measurable set with $m(E)\ge \lambda>0$, with $m$ being the Lebesgue measure. Then
\begin{align*}
\|T_d\|_{C([0,2\pi])} \le \frac{1}{2}\pth{\frac{2}{\sin(\lambda/4)}}^{2d}\|T_d\|_{C(E)}. 
\end{align*}
\end{lemma}

Choosing $d=2$ in \Cref{lemma:remez} and using that $\|A\|_{C[0,2\pi]}\ge \frac{2nr^2}{(1+r^2)^2}$, we get
\begin{align*}
m\pth{\sth{ \theta\in [0,2\pi]: A(\theta)\le a\cdot \frac{2nr^2}{(1+r^2)^2}}} \le 4\arcsin\pth{(8a)^{1/4}} \le 2\pi(8a)^{1/4}
\end{align*}
for all $a\in [0,\frac{1}{8}]$, where the second inequality uses $\arcsin x\le \frac{\pi}{2}x$ for $x\ge 0$. The above inequality trivially extends to $a>\frac{1}{8}$ as well, so it holds for all $a\ge 0$. Therefore, for every $t\ge 0$, 
\begin{align}\label{eq:sublevel}
m\pth{\sth{\theta\in [0,2\pi]: B(\theta)\le t }} &\le m\pth{\sth{\theta\in [0,2\pi]: A(\theta)\le t(\sqrt{t}+2)^2 }} \nonumber \\
&\le 2\pi\cdot \pth{\frac{8t(\sqrt{t}+2)^2}{\frac{2nr^2}{(1+r^2)^2}}}^{1/4}. 
\end{align}
By \eqref{eq:q-upper-bound} and \eqref{eq:sublevel}, we integrate the tail to get
\begin{align*}
\frac{1}{2\pi}\int_0^{2\pi} q(\theta)\rmd \theta &\le \frac{1}{2\pi}\int_0^{2\pi} \exp\pth{-\frac{B(\theta)}{2}}\rmd \theta \\
&= \frac{1}{4\pi}\int_0^\infty e^{-t/2}m\pth{\sth{\theta\in [0,2\pi]: B(\theta)\le t }} \rmd t \\
&\le \frac{1}{2}\pth{\frac{(1+r^2)^2}{2nr^2}}^{1/4}\int_0^\infty e^{-t/2}(8t(\sqrt{t}+2)^2)^{1/4}\rmd t \\
&\le C_1\pth{\frac{(1+r^2)^2}{nr^2}}^{1/4}, 
\end{align*}
with 
\begin{align*}
C_1 = \frac{1}{2^{5/4}}\int_0^\infty e^{-t/2}(8t(\sqrt{t}+2)^2)^{1/4}\rmd t < 2.812.
\end{align*}
Finally, by \eqref{eq:basic-ineq} and Stirling's approximation $\sqrt{2\pi n}(\frac{n}{e})^n\le n!\le e^{1/12}\sqrt{2\pi n}(\frac{n}{e})^n$ for every $n\ge 1$, we plug in the expression of $r=\sqrt{\frac{k}{n-k}}$ and get the desired result
\begin{align*}
|e_k(z)| &\le C_1 \frac{(1+r^2)^{(n+1)/2}}{n^{1/4}r^{k+1/2}} \\
&= C_1\pth{\frac{n^{n+1/2}}{k^{k+1/2}(n-k)^{n-k+1/2}}}^{1/2} \\
&\le C_1\pth{ \sqrt{2\pi}e^{1/6} \binom{n}{k} }^{1/2} \le \sqrt{C\binom{n}{k}}, 
\end{align*}
with $C = C_1^2\sqrt{2\pi} e^{1/6} < 24$. 

\subsection{Proof of \Cref{thm:ESP-matrix}}
Since $e_1(A)=0$, we assume that $k\ge 2$. The proof uses a generating polynomial that respects both row and column sum conditions of $A$. Let $p\in [0,1]$, $X=(X_1,\dots,X_n)$ be a random vector in $\bR^n$ with i.i.d. $\Bern(p)$ components, and $Y$ be an independent copy of $X$. Consider
\begin{align*}
Q_k(p) &:= \bE\qth{(X^\top AY)^k} = \sum_{i_1,j_1,\dots,i_k,j_k} \prod_{\ell=1}^k A_{i_\ell,j_\ell} \bE\qth{ \prod_{\ell=1}^k X_{i_\ell} }\bE\qth{ \prod_{\ell=1}^k Y_{j_\ell} }. 
\end{align*}
Since $\bE[X_i^r] = \bE[Y_j^r] =p$ for all $r\ge 1$, $Q_k(p)$ is a polynomial in $p$ of degree at most $2k$, with leading coefficient given by
\begin{align*}
\sum_{\substack{i_1,\dots,i_k \text{ distinct} \\ j_1,\dots,j_k \text{ distinct}}} \prod_{\ell=1}^k A_{i_\ell,j_\ell} = (k!)^2 e_k(A).
\end{align*}
In addition, since $A{\bf 1}=A^\top {\bf 1}=0$, we have
\begin{align*}
Q_k(p) = \bE\qth{\pth{(X-p{\bf 1})^\top A (Y-p{\bf 1}) }^k} =: \bE\qth{ (\xi^\top A \zeta)^k }.
\end{align*}
Our goal is to provide a uniform upper bound of $|Q_k(p)|$ on $[0,1]$, and use it to upper bound the leading coefficient $(k!)^2 e_k(A)$. 

First, we condition on $\zeta$ and take expectation with respect to the centered vector $\xi$. By Hoeffding's lemma, $\xi$ is $\frac{1}{4}$-subGaussian, so $\xi^\top A\zeta$ is $\frac{\|A\zeta\|_2^2}{4}$-subGaussian. We use the following moment bound for subGaussian random variables: by standard tail integration, 
\begin{align}\label{eq:moment}
\bE\qth{|X|^k} \le 2k\int_0^\infty x^{k-1}e^{-2x^2/\sigma^2} \rmd x = 2^{k/2+1}\sigma^k\Gamma(\frac{k}{2}+1) \le 2(k\sigma^2)^{k/2}, 
\end{align}
where the second inequality uses $\Gamma(\frac{k}{2}+1)\le (\frac{k}{2})^{k/2}$ for $k\ge 2$. By \eqref{eq:moment}, 
\begin{align*}
\bE \qth{ |\xi^\top A\zeta|^k } \le 2\pth{ \frac{k}{4} }^{k/2} \bE\qth{\|A\zeta\|_2^k}. 
\end{align*}
To proceed, note that for each coordinate $i\in [n]$, the map $\{0,1\}\ni \zeta_i\mapsto \|A\zeta\|_2$ has the bounded difference property with difference at most $\|A_{i,:}\|_2$, where $A_{i,:}$ is the $i$-th row vector of $A$. Since $\sum_{i=1}^n \|A_{i,:}\|_2^2 = \|A\|_{\mathrm{F}}^2 = n^2$, McDiarmid's inequality gives
\begin{align*}
\bP\pth{ \left| \|A\zeta\|_2 - \bE[\|A\zeta\|_2] \right| \ge t } \le \exp\pth{-\frac{2t^2}{n^2}} \quad \text{for all }t\ge 0. 
\end{align*}
Therefore, by the same tail integration in \eqref{eq:moment}, 
\begin{align*}
\bE\qth{ \left| \|A\zeta\|_2 - \bE[\|A\zeta\|_2] \right|^k } \le 2(k\frac{n^2}{4})^{k/2}  = 2\pth{ \frac{k}{4} }^{k/2} n^k. 
\end{align*}
In addition, 
\begin{align*}
\bE[\|A\zeta\|_2] \le \sqrt{\bE\|A\zeta\|_2^2} = \sqrt{\|A\|_{\mathrm{F}}^2 p(1-p)} \le \frac{n}{2}. 
\end{align*}
Consequently, by triangle inequality, 
\begin{align*}
\bE\qth{\|A\zeta\|_2^k} \le 2^{k-1}\pth{ \bE\qth{ \left| \|A\zeta\|_2 - \bE[\|A\zeta\|_2] \right|^k } + \bE[\|A\zeta\|_2]^k } \le \pth{k^{k/2} + \frac{1}{2}}n^k \le 2k^{k/2}n^k. 
\end{align*}
In summary, we have shown that for all $p\in [0,1]$, 
\begin{align}\label{eq:norm-bound}
|Q_k(p)| \le 4\pth{\frac{k}{2}}^k n^k. 
\end{align}

Finally, we recall the well-known result for Chebyshev polynomials \cite[Chapter 3, Theorem 6.1]{devore1993constructive}: if a polynomial $P$ of degree $m$ satisfies $\|P\|_{C[-1,1]}\le 1$, the largest possible leading coefficient is attained by the degree-$m$ Chebyshev polynomial, and this value is $2^{m-1}$. Applying this result to \eqref{eq:norm-bound} with a linear transformation $x\mapsto \frac{x+1}{2}$ from $[-1,1]$ to $[0,1]$, then
\begin{align*}
\frac{(k!)^2 |e_k(A)|}{2^{2k}} \le 2^{2k-1}\cdot 4\pth{\frac{k}{2}}^k n^k, 
\end{align*}
which simplifies to
\begin{align*}
|e_k(A)| \le \frac{2 (8nk)^k}{(k!)^2} \le \frac{(8e)^k n^k}{k!} \le (8e^2)^k \binom{n}{k}, 
\end{align*}
by Stirling's approximation $k!\ge \sqrt{2\pi k}\pth{\frac{k}{e}}^k$ and the inequality $\frac{n^k}{n(n-1)\cdots(n-k+1)}\le e^k$. This is the desired result.

\subsection{Proof of \Cref{cor:permanent}}
For a PSD and doubly stochastic matrix $M$, the centered matrix $A:=M-\frac{J}{n}$ remains PSD, where $J$ is the all-ones matrix, and $A$ has zero row and column sums. By \cite[Appendix B]{han2024approximate}, it holds that
\begin{align*}
\frac{n^n}{n!}\Perm(M) = \sum_{k=0}^n \frac{n^k}{\binom{n}{k}} e_k(A),  
\end{align*}
where $e_k(A)$ is given in \eqref{eq:ESP-matrix}. We prove the two upper bounds separately: 
\begin{enumerate}
    \item For the first upper bound, recall from \cite[Lemma 5.7, Lemma B.1]{han2024approximate} that
    \begin{align*}
    e_k(A) = \frac{1}{k!} \bE\qth{|e_k(z)|^2}, \quad z\sim \calC\calN(0,A). 
    \end{align*}
    Here $\calC\calN(0,A)$ denotes a complex normal random vector with mean zero and covariance $A$. By the eigendecomposition of $A=UDU^\top$ with $D=\diag(0,\lambda_2,\dots,\lambda_n)$, it is clear that $z$ could be written as $UD^{1/2}w$, with $w=(w_1,\dots,w_{n})\sim \calC\calN(0,I_{n})$. In particular, since $U$ is unitary, we have $\|z\|_2^2 = \|D^{1/2}w\|_2^2 =  \sum_{i=2}^n \lambda_i |w_i|^2$. By \Cref{thm:ESP-vector}, 
    \begin{align*}
    \bE\qth{|e_k(z)|^2} &\le C\binom{n}{k} \bE\bqth{\pth{\frac{\|z\|_2^2}{n}}^k} = C\binom{n}{k} \bE\bqth{\pth{\frac{\sum_{i=2}^n \lambda_i|w_i|^2}{n}}^k} \\
    &= C\frac{k!}
    {n^k}\binom{n}{k} \sum_{ \substack{k_2,\dots,k_n\ge 0 \\ k_2+\dots+k_n=k }} \lambda_2^{k_2}\cdots\lambda_n^{k_n}\bE\bqth{\frac{|w_2|^{2k_2}\cdots |w_n|^{2k_n}}{k_2! \cdots k_n!}} \\
    &= C\frac{k!}
    {n^k}\binom{n}{k} \sum_{ \substack{k_2,\dots,k_n\ge 0 \\ k_2+\dots+k_n=k }} \lambda_2^{k_2}\cdots\lambda_n^{k_n}, 
    \end{align*}
    where the last step follows from $\bE\qth{|w|^{2k}} = k!$ for $w\sim \calC\calN(0,1)$. Rearranging gives the first upper bound. 
    \item For the second upper bound, \Cref{thm:ESP-matrix} gives
    \begin{align*}
    |e_k(A)| \le \binom{n}{k} \pth{\frac{B\|A\|_{\mathrm{F}}}{n}}^k = \binom{n}{k} \bpth{\frac{B\sqrt{\lambda_2^2+\cdots+\lambda_n^2}}{n}}^k.
    \end{align*}
    The claimed upper bound then follows from simple algebra.
\end{enumerate}

\subsection{Proof of \Cref{cor:chi-squared}}
Let $M\in \bR^{n\times n}$ be the channel overlap matrix generated by $P_1,\dots,P_n$, \cite[Lemma 5.1]{han2024approximate} gives
\begin{align*}
\chi^2(\bP_n \| \bQ_n) = \frac{n^n}{n!}\Perm(M) - 1. 
\end{align*}
In addition, by \cite[Lemma 5.2]{han2024approximate}, the matrix $M$ is PSD and doubly stochastic. Therefore, by the first inequality of \Cref{cor:permanent}, 
\begin{align*}
\chi^2(\bP_n \| \bQ_n) &\le C\sum_{k=2}^n \sum_{ \substack{k_2,\dots,k_n\ge 0 \\ k_2+\dots+k_n=k }} \lambda_2(M)^{k_2}\cdots\lambda_n(M)^{k_n} \\
&\le C\sum_{ \substack{k_2,\dots,k_n\ge 0 \\ k_2+\dots+k_n\ge 2 }} \lambda_2(M)^{k_2}\cdots\lambda_n(M)^{k_n} \\
&= C\bpth{\sum_{k_2,\dots,k_n\ge 0} \lambda_2(M)^{k_2}\cdots\lambda_n(M)^{k_n}  - 1 - \sum_{i=2}^n \lambda_i(M)} \\
&= C\bpth{ \prod_{i=2}^n \frac{1}{1-\lambda_i(M)} - 1 - \sum_{i=2}^n \lambda_i(M) }. 
\end{align*}

\subsection{Proof of \Cref{cor:deFinetti}}\label{subsec:deFinetti}
Let $A:=M-\frac{J}{n}$, where $M$ is the channel overlap matrix and $J$ is the all-ones matrix. By \cite[Appendix D.3]{han2024approximate}, it holds that
\begin{align*}
\chi^2(\bP_{n,k}\|\bQ_{n,k}) = \sum_{\ell=2}^k \frac{\binom{k}{\ell}}{\binom{n}{\ell}} \frac{n^{\ell}}{\binom{n}{\ell}}e_\ell(A) &\le \sum_{\ell=2}^k \pth{\frac{k}{n}}^{\ell} \frac{n^{\ell}}{\binom{n}{\ell}}e_\ell(A) \\
&= \sum_{\ell=2}^k \frac{n^{\ell}}{\binom{n}{\ell}}e_\ell\pth{\frac{k}{n}A} 
= \frac{n^n}{n!}\Perm\pth{\frac{k}{n}A + \frac{J}{n}}-1. 
\end{align*}
Since $M' = \frac{k}{n}A+\frac{J}{n}$ is still a doubly stochastic and PSD matrix, applying \Cref{cor:permanent} to $M'$ gives
\begin{align*}
\frac{n^n}{n!}\Perm\pth{\frac{k}{n}A + \frac{J}{n}}-1 &\le \sum_{\ell=2}^n B^{\ell} \pth{\lambda_2^2(M') + \cdots + \lambda_n^2(M')}^{\ell/2} \\
&= \sum_{\ell=2}^n \pth{\frac{Bk}{n}\sqrt{\lambda_2^2(M) + \cdots + \lambda_n^2(M)}}^{\ell}.
\end{align*}
To proceed, note that $\sum_{i=2}^n \lambda_i^2(M)\le \pth{\sum_{i=2}^n \lambda_i(M)}^2 \le {\sf T}^2 $, and $\sum_{i=2}^n \lambda_i^2(M)\le \sum_{i=2}^n \lambda_i(M) = {\sf T}$. If $\frac{k}{n}\sqrt{{\sf T}^2 \wedge {\sf T}} \le \frac{1}{2B}$, the above geometric series is dominated by the contribution of $\ell=2$, which gives the claimed result.

We also prove that ${\sf T}\le |\calX|-1$ when $P_1,\dots,P_n$ are supported on a finite set $\calX$. The channel overlap matrix $M$ is defined as
\begin{align}\label{eq:channel-overlap}
M_{ij} = \frac{1}{n}\int \frac{\rmd P_i \rmd P_j}{\rmd \overline{P}} = \sum_{x\in \calX} \frac{P_i(x)P_j(x)}{\sum_{\ell=1}^n P_{\ell}(x)}. 
\end{align}
Therefore, 
\begin{align*}
\mathrm{Tr}(M) = \sum_{i=1}^n M_{ii} = \sum_{x\in \calX} \frac{\sum_{\ell=1}^n P_{\ell}(x)^2}{\sum_{\ell=1}^n P_{\ell}(x)} \le \sum_{x\in \calX} 1 = |\calX|,
\end{align*}
and ${\sf T} = \mathrm{Tr}(M) - 1 \le |\calX| - 1$. 

\bibliographystyle{alpha}
\bibliography{refs}

\end{document}